\documentclass[12pt]{article}
\usepackage{mathrsfs}
\usepackage{amsfonts}
\usepackage{amssymb,amsfonts,amsmath, psfrag,eepic}
\newcommand{\D}{\displaystyle}

\@addtoreset{equation}{section} \@addtoreset{theo}{section}
\makeatother

\begin{document}
\begin{center}
{\large \bf The Riordan Group and Symmetric Lattice Paths}
\end{center}

\begin{center}
Li-Hua Deng$^1$, Eva Y. P. Deng$^2$ and Louis W. Shapiro$^3$
\end{center}
\vspace{0.2cm} \centerline{$^{1,\  2}$School of Mathematical Sciences, Dalian University of Technology}
\centerline{Dalian, 116024, P. R. China}
\vspace{3mm} \centerline{$^{3}$Department of Mathematics, Howard University,}
\centerline{Washington, DC 20059, USA}
 \vspace{3mm}
\centerline{$^1$denglihua840321@163.com, $^2$ypdeng@dlut.edu.cn,
$^3$lshapiro@howard.edu}

\vspace{0.5cm}

\noindent \textbf{Abstract:}\ \ In this paper, we study symmetric
lattice paths. Let $d_{n}$, $m_{n}$, and $s_{n}$ denote the number
of symmetric Dyck paths, symmetric Motzkin paths, and symmetric Schr\"{o}der
paths of length $%
2n$, respectively. By using Riordan group methods we obtain six
identities relating $d_{n}$, $m_{n}$, and $s_{n}$ and also give two
of them combinatorial proofs. Finally, we investigate some relations
satisfied by the generic element of some special Riordan arrays and get the average mid-height and the
average number of points on the $x$-axis of symmetric Dyck paths of
length $2n.$

\noindent \textbf{Key words:} Symmetric Dyck paths; Symmetric
Motzkin paths; Symmetric Schr\"{o}der paths; Riordan group;
Combinatorial identities

\noindent{\textbf{1\hspace{0.4cm} Introduction} }

Lattice paths have been widely studied from various points of view.
A surprisingly large number of $1-1$ correspondences are known that
relate these lattice paths to other classes of objects, such as
trees, polygon triangulations, 213 avoiding permutations, $2\times
n$ standard Young tableaux, and so on. A large number of references
can be found in [1,2,12]. In the present paper we will study some
symmetric lattice paths.

Let $k$ be any fixed positive integer. In the plane
$\mathbb{Z}\times \mathbb{Z}$, we consider lattice paths with three
step types: \emph{an up step} $U=(1,1)$, \emph{a down step}
$D=(1,-1)$ and a \emph{$k$-horizontal step} $(k,0)$. Usually a
$1$-horizontal step is simply called \emph{a horizontal step}. For
convenience, we denote $h=(1,0),$ and $H=(2,0).$

A \emph{generalized Motzkin path} of length $n$ is a lattice path
from
origin $(0,0)$ to $(n,0)$ consisting of up steps, down steps and $k$%
-horizontal steps that never goes below the $x$-axis. In the general
setting of this paper, when $k=0,1,2$, we call it a Dyck path, a
Motzkin path and a Schr\"{o}der path, respectively. Recall that the
generating functions for the number of Dyck paths, Motzkin paths
and Schr\"{o}der paths are
$\D C\left( x\right) =\frac{1-\sqrt{1-4x}}{2x},M\left( x\right) =\frac{1-x-\sqrt{%
1-2x-3x^{2}}}{2x^{2}}\
 \text{and} \ R\left( x\right) =\D\frac{1-x-\sqrt{%
1-6x+x^{2}}}{2x}$.

Another generating function that will occur often is
\begin{equation*}
B\left( x\right) =\frac{1}{\sqrt{1-4x}}=\sum_{n\geq 0}\binom{2n}{n}%
x^{n}=1+2x+6x^{2}+20x^{3}+70x^{4}+\cdots .
\end{equation*}%
Since $B\left( x^{2}\right) $ is the generating function for paths from $%
\left( 0,0\right) $ and ending on the $x$-axis using just up and
down steps, we see that
\begin{equation*}
B\left( x^{2}\right) =\frac{1}{1-2x^{2}C\left( x^{2}\right)
}=\frac{C\left( x^{2}\right) }{1-x^{2}C^{2}\left( x^{2}\right) }.
\end{equation*}%
For the first equality the 2 comes since after each return to the
$x$-axis the next step can be either up or down. The second equality
follows from starting a new component each time the path goes from
the $x$-axis down to the line $y=-1.$

A closely related generating function identity, perhaps of
independent interest, is
\begin{equation*}
\frac{1}{1-2x}=\frac{C\left( x^{2}\right) }{\left( 1-xC\left(
x^{2}\right) \right) ^{2}}.
\end{equation*}

A \emph{generalized symmetric Motzkin path} of length $2n$ is a
lattice path such that if the $i$-th step of generalized Motzkin
path is an up step (down step or $k$-horizontal step) then the
$2n-i$-th step is down step (up step or $k$-horizontal step).
Basically the path is the same when read from left to right as when
read from right to left with a vertex in the middle. In this paper,
when $k=0,1,2$, we will call it a symmetric Dyck path, a symmetric
Motzkin path and a symmetric Schr\"{o}der path, respectively.

The number of symmetric Dyck paths of length $2n$ equals the $n$-th
Central binomial coefficient (one version, this one is sequence
A001405 in [8]), denoted by $d_{n}$ . The first few numbers are
$1,1,2,3,6,10,20,35,\cdots $ and the generating function of
$({d}_{n})_{n\in \mathbb{N}}$ is
\begin{equation*}
{d}(x)=\displaystyle\frac{1}{2x}\left(
\sqrt{\frac{1+2x}{1-2x}}-1\right) =B\left( x^{2}\right) +xB\left(
x^{2}\right) C\left(
x^{2}\right)=\frac{C\left(x^2\right)}{1-xC(x^2)} .
\end{equation*}

The number of symmetric Motzkin paths of length $2n$ equals the
number of directed animals of size $n+1$ (sequence A005773 in [8]),
denoted by $m_{n}$. The first few numbers are 1, 2, 5, 13, 35, 96, 267, 750 $\dots$ and the generating function of $({m}_{n})_{n\in
\mathbb{N}}$ is
\begin{equation*}
{m}(x)=\displaystyle\frac {1}{2x}\left(\sqrt{\frac
{1+x}{1-3x}}-1\right).
\end{equation*}

Let $s_{n}$ denote the number of symmetric Schr\"{o}der paths of length $2n$. An equation for the generating function $s(x)=\sum_{n\geq
0}s_{n}x^{n}$ is obtained from the \textquotedblleft first return
decomposition" of a symmetric Schr\"{o}der path $P$. If the first
step of $P$ is a 2-horizontal step, then the contribution of this
case gives $x^{2}s(x)$. If the first
step of $P$ is an up step, then the contribution of this case gives $
xs(x)+x^{2}s(x)R(x^{2})$, where
$R(x)$ is the generating function
for Schr\"{o}der paths. The term $xs\left( x\right) $ accounts for
those paths that do not touch the $x$-axis until the last step while
the term $x^{2}s\left( x\right) R\left( x^{2}\right) $ accounts for
paths that do return, $xR\left( x^2\right) x$ gives the first return and then $
s\left( x\right) $ for the part between the first and last returns.
Hence
\begin{equation*}
s(x)=1+xs(x)+x^{2}s(x)+x^{2}s(x)R(x^{2})
\end{equation*}%
which implies that
\begin{equation*}
\displaystyle{s}(x)=\displaystyle\frac{1}{2x}\left( \sqrt{\frac{x^{2}-2x-1}{%
x^{2}+2x-1}}-1\right) .
\end{equation*}%
In fact, $s_{n}$ is the sequence A026003 in [8]. The first few
numbers are
1, 1, 3, 5, 13, 25, 63, 129. These three sequences are all mentioned in $%
\left[ 8\right] $ and were contributed in 2002-3 by Emeric Deutsch
although without proofs. An equivalent way to view these sequences
is as left factors of Dyck (Motzkin, Schr\"{o}der) paths. By taking
the mirror image of the left factor as the right factor you obtain
the symmetric path.

Another natural way to define Symmetric paths would be as a path
that is symmetric about the line, $x=n$ or as a path that looks the
same whether
going from left to right or from right to left. Call such paths \textbf{%
palindromic}. There exist palindromic Motzkin paths of odd length
and palindromic Schr\"{o}der paths with a 2-horizontal step from $x=n-1$ to
$\ n+1.$
Our definition however requires a lattice point on the path when \thinspace $%
x=n.$ Palindromic Motzkin paths of odd length, $2n+1,$ must have a
horizontal step in the middle and removing it gives a symmetric Motzkin
paths of length $2n.$ Similarly palindromic Schr\"{o}der paths with
a 2-horizontal middle step can map to a symmetric Schr\"{o}der path by
removing the middle step. Thus our results are easily translated
over to the palindromic version.

The sequences above occur in other combinatorial structures such as
symmetric ordered trees. However in the this paper we describe them
in terms of symmetric lattice paths. In Section 3 we will give six
identities related to them by using Riordan group
techniques which will be reviewed in Section 2. In Section 4 we will
study some relations satisfied by the generic element of some
special Riordan arrays and get the average mid-height and the average
number of points on the $x$-axis of symmetric Dyck paths.

\noindent{\textbf{2\hspace{0.4cm} Riordan Group } }{\ }

In 1978, Rogers [6] introduced the renewal array, which is a
generalization of the Pascal, and Motzkin triangles. Kettle [5] used
the theory of renewal arrays to study other types of combinatorial
triangles, especially those found in walk problems. Shapiro et al.
[9] and Sprugnoli [11] generalized these kind of arrays to Riordan
arrays and the Riordan group. Riordan arrays constitute a practical
device for solving combinatorial sums by means of composition of
generating functions.

\emph{A Riordan array }is an infinite lower triangular array $R=\{{r_{n,k}}%
\}_{n\geq k\geq 0}$ generated by a pair of analytic functions $%
g(x)=1+g_{1}x+g_{2}x^{2}+\cdots $ and $f(x)=f_{1}x+f_{2}x^{2}+\cdots
.$ If also $f_{1}\neq 0$ then we have an element of the Riordan
group. The array is defined by
\begin{equation*}
r_{n,k}=[x^{n}]g(x)f(x)^{k},
\end{equation*}%
where the notion $[x^{n}]$ denotes the \textquotedblleft coefficient
operator" that extracts the coefficient of $x^{n}$. We often denote
a Riordan array as $R=(g(x),f(x))$ or even as $(g,f).$

Suppose we multiply the matrix $R=(g,f)$ by a column vector $%
(a_{0},a_{1},\cdots )^{T}$ and get a column vector
$(b_{0},b_{1},\cdots
)^{T} $. Let $A(x)$ and $B(x)$ be the generating functions for the sequence $%
(a_{0},a_{1},\cdots )$ and $(b_{0},b_{1},\cdots )$ respectively.
Then it follows quickly that
\begin{equation*}
B(x)=(g(x),f(x))\ast A(x)=g(x)A(f(x)).
\end{equation*}%
This is the essential fact sometimes referred to as ``The Fundamental
Theorem of Riordan Arrays" or even as the FTRA. Many examples and
properties of the Riordan group are described in [9,10,11] along
with the connection to the Lagrange inversion formula.

The \emph{Riordan group } $\mathcal{R}$ = \{${R|R=(g(x),f(x))}$
consists of the Riordan arrays with $g_{0}=1,f_{0}=0,f_{1}\neq 0$
\}. The multiplication in $\mathcal{R}$ is just matrix
multiplication and is given by
\begin{equation*}
(g_{1}(x),f_{1}(x))\ast
(g_{2}(x),f_{2}(x))=(g_{1}(x)g_{2}(f_{1}(x)),f_{2}(f_{1}(x))).
\end{equation*}%
The identity is $I=(1,x)$, and the inverse matrix is specified by
\begin{equation*}
(g(x),f(x))^{-1}=\left( \frac{1}{g(\bar{f}(x))},\bar{f}(x)\right) ,
\end{equation*}%
where $\bar{f}$ is the compositional inverse of $f(x).$ That is
\begin{equation*}
f(\bar{f}(x))=\bar{f}(f(x))=x.
\end{equation*}%
\quad It is easy to see that the Riordan group can used to study
inverse relations. From $R$ and $R^{-1}$, an inverse relation can be
established and thus we have a systematic way to find inverse
relations and sums. The paper [9] and [10] provide many examples. It
should be noticed that earlier Riordan [7] and Gould and Hsu [4]
studied many inverse relations.

\noindent{\textbf{3\hspace{0.4cm} Combinatorial Identities} }{\ }

Deng and Yan [3] obtained some identities involving
the Catalan, Motzkin and Schr\"{o}der numbers using Riordan group
methods. In this section we will give six identities among $d_{n}$,
$m_{n}$, and $s_{n}$ by using the Riordan group and give
combinatorial proofs for two of them. We start by noting a very
suggestive equation connecting the first few terms
of $m_{n}$ and $d_{n}.$%
\begin{equation*}
\left[
\begin{array}{cccccc}
1 & &  &  & &  \\
1 & 1 & \\
1 & 2 & 1 &   \\
1 & 3 & 3 & 1 &  \\
1 & 4 & 6 & 4 & 1 &  \\
  &  & \cdots &  &  & \ddots%
\end{array}%
\right] \left[
\begin{array}{c}
1 \\
1 \\
2 \\
3 \\
6 \\
\vdots
\end{array}%
\right] =\left[
\begin{array}{c}
1 \\
2 \\
5 \\
13 \\
35 \\
\vdots%
\end{array}%
\right]
\end{equation*}%
Since the first matrix above is the Pascal matrix
 $$P=\left(\frac{1}{1-x},\frac{x}{1-x}\right),$$
we can use the fundamental theorem as follows.

\noindent {\textbf{Theorem~3.1.}}\ For $n\geq 0$, we have the
following formula
\begin{equation}
\displaystyle {m}_{n}=\sum_{k=0}^{n} {\binom{n}{k}}{d}_{k}.
\end{equation}

\noindent \textbf{Proof.}\ Since the generating function of sequence $({d}%
_{n})_{n\in \mathbb{N}}$ is ${d}(x),$ it follows from the
fundamental theorem that
\begin{equation*}
\displaystyle\left(\frac{1}{1-x},\frac{x}{1-x}\right)*
d(x)=\frac{1}{1-x}{d}\left( \frac{x}{1-x}\right) ={m}(x).
\end{equation*}%
So we get
\begin{equation*}
\displaystyle{m}_{n}=\sum_{k=0}^{n}{\binom{n}{k}}d_{k}.
\end{equation*}

\noindent \textbf{Combinatorial proof.}\ Suppose that a symmetric
Motzkin path of length $2n$ contains $2k$ horizontal steps. We can
reduce it to a symmetric Dyck path of length $2(n-k)$ by removing
all the horizontal steps. Conversely, given a symmetric Dyck path of
length $2(n-k)$, we can reconstruct $\binom{n}{k}$ symmetric Motzkin
paths of length $2n$ by inserting $2k$ horizontal steps. The
symmetric condition means placing $k$ horizontal steps into the left
half of the Dyck path. There are $n-k+1$ vertices in which we can
place $k$ horizontal steps with repetition allowed so we have
$\binom{\left( n-k+1\right) +k-1}{k}=\binom{n}{k}$ possibilities.
Therefore we do have the relation
\begin{equation*}
\displaystyle{m}_{n}=\sum_{k=0}^{n}{\binom{n}{k}}d_{k}.\quad
_{\blacksquare }
\end{equation*}

\noindent Now we get the inverse of the identity (1) by multiplying
the Riordan matrix inverse
\begin{equation*}
P^{-1}=\left( \frac{1}{\frac{1}{1-\left( \frac{x}{1+x}\right) }},\frac{\frac{%
x}{1+x}}{1-\left( \frac{x}{1+x}\right) }\right)=\left( \frac{1}{1+x},\frac{x%
}{1+x}\right) =\left[
\begin{array}{cccccc}
1 &  &  &  &  \\
-1 & 1 &  &  &  \\
1 & -2 & 1 &  &  \\
-1 & 3 & -3 & 1 &  \\
1 & -4 & 6 & -4 & 1 &\\
&  &\cdots &  &  & \ddots%
\end{array}%
\right] .
\end{equation*}%

\noindent So we have

 \noindent{\textbf{Theorem~3.2.}}\ For $n\geq 0$, we have the following formula
\begin{equation}
\displaystyle {d}_{n}=\sum_{k=0}^{n}(-1)^{n-k}
{\binom{n}{k}}{m}_{k}.
\end{equation}

\noindent {\textbf{Remark:}}\ In fact, the formulas (1) and (2) are a
special case of the inverse transformation

\begin{equation*}
\quad\displaystyle b_{n}=\sum_{k=0}^{n} {\binom{n}{k}}a_{k}
\Longleftrightarrow \displaystyle a_{n}=\sum_{k=0}^{n}(-1)^{n-k} {\binom{n}{k%
}}b_{k}.
\end{equation*}

The transformation is widely used in the study of integer sequences
where it is called the binomial transform. It is also called the
Euler transform as it was introduced by Euler as a tool to
accelerate the speed of convergence of sequences.

\noindent {\textbf{Theorem~3.3.}}\ For $n\geq 0$, we have the
following formula
\begin{equation}
\displaystyle {s}_{n}=\sum_{k=0}^{\lfloor n/2\rfloor }
{\binom{n-k}{n-2k}}{d}_{n-2k}.
\end{equation}

\noindent \textbf{Proof.}\ Consider the Riordan array

\begin{equation*}
D=\displaystyle\left( \frac{1}{1-x^{2}},\frac{x}{1-x^{2}}\right)
=\left[
\begin{array}{ccccccc}
1 &  &  &  &  &  &  \\
0 & 1 &  &  &  &  &  \\
1 & 0 & 1 &  &  &  &  \\
0 & 2 & 0 & 1 &  &  &  \\
1 & 0 & 3 & 0 & 1 &  &  \\
0 & 3 & 0 & 4 & 0 & 1 &  \\
&  &  \ldots &  &  &  & \ddots%
\end{array}%
\right]
\end{equation*}%

\noindent Since the generic element of the Riordan array $ \D \left(
\frac{1}{1-x^{2}},\frac{x}{1-x^{2}}\right)$ is

$$\displaystyle[x^{n}]\frac{x^{j}}{(1-x^{2})^{j+1}}=\left\{
\begin{array}{cc}
{\binom{\frac{n+j}{2} }{\frac{n-j}{2}}} & n-j\text{ even} \\
0 & n-j\text{ odd}%
\end{array}%
\right. ,$$
we finish the proof by setting
$\frac{n-j}{2}=k$.
\bigskip

\noindent \textbf{Combinatorial proof.}\ Suppose that a symmetric Schr\"{o}%
der path ending at (2n,0) contains $2k$ 2-horizontal steps.  We can reduce it to a symmetric Dyck path of length $2(n-2k)$
by deleting all 2-horizontal steps. Conversely, given a symmetric
Dyck path of length $2(n-2k)$, we can construct ${\binom{n-k}{k}}$
symmetric Schr\"{o}der paths of length $2n$ by inserting $2k$
2-horizontal steps. So the numbers ${s}_{n}$ relate to the numbers
${d}_{n}$ by
\begin{equation*}
\displaystyle{s}_{n}=\sum_{k=0}^{\lfloor n/2\rfloor }{n-k\choose
n-2k}d_{n-2k}.\quad _{\blacksquare }
\end{equation*}

\noindent {\textbf{Theorem~3.4.}}\ For $n\geq 0$, we have the
following formula
  \begin{equation}
  {d}_{n}=\sum_{k=0}^{n}(-1)^{\left( n-k\right) /2}d_{n,k}{s}_{k},
  \end{equation}%
where
  $d_{n,k}=\left\{\begin{array}{cc}
  \frac{k+1}{n+1}{\binom{n+1}{\frac{n-k}{2}}} & n-k\text{ even} \\
  0 & n-k\text{ odd}%
  \end{array}\right. ,$
and the generating function for
  $({d}_{n,k})_{n\in \mathbb{N}}$
is
  $x^k C^{k+1}(x^2),
  k= 0, 1,2,\cdots.$

\noindent \textbf{Proof.}\ Since
 $$\begin{array}{ll} \displaystyle
 D^{-1}&=\displaystyle\left( \frac{1}{1-x^{2}},\frac{x}{1-x^{2}}\right)^{-1}=\left( \displaystyle\frac{\sqrt{1+4x^{2}}-1}{2x^{2}},%
 \frac{\sqrt{1+4x^{2}}-1}{2x}\right)
 \end{array}$$

\noindent The generic element of $D^{-1}$ is
\begin{equation*}
\begin{array}{ll}
& [x^{n}]x^{k}\displaystyle\left(
\frac{\sqrt{1+4x^{2}}-1}{2x^{2}}\right)
^{k+1} =\left\{
\begin{array}{ll}
\displaystyle(-1)^{\left( n-k\right) /2}\frac{k+1}{n+1}{\binom{n+1}{\frac{n-k%
}{2}}}\ \ \ \ \ \ n-k\ \ \mathrm{even}; \\
&  \\
\quad \quad \quad \quad 0\ \ \ \ \ \ \ \ \quad \quad \quad \quad
\quad \quad \ \quad \mathrm{otherwise}. &
\end{array}%
\right.%
\end{array}%
\end{equation*}%
Set  $d_{n,k}=\frac{k+1}{n+1}{\binom{n+1}{\frac{n-k}{2}}}$,
then we obtain the formula (4).
It is easy to see that $d_{n,k}$ is the $\left( n,k\right)
$ entry of the Riordan array
\begin{equation*}
D^{\ast }=\left( \displaystyle\frac{1-\sqrt{1-4x^{2}}}{2x^{2}},\frac{1-\sqrt{%
1-4x^{2}}}{2x}\right) =\left(C(x^2),xC(x^2)\right) \quad _{\blacksquare }
\end{equation*}

\noindent{\textbf{Theorem~3.5.}}\ For $n\geq 0$, we have the following
formula
 \begin{align}
  \quad {s}_{n}=\sum_{k=0}^{n}(-1)^{n-k} {s}_{n,k}{m}_{k},
 \end{align}
where
 $s_{n,k}= \sum_{j=0}^{n-k}{\binom{k+j}{k}}{\binom{j}{n-j-k}},$
and the generating function of
 $({s}_{n,k})_{n\in \mathbb{N}}$
is\\
 $x^{k}\left( \frac{1}{1-x-x^{2}}\right) ^{k+1}$, $k=0,1,2,\cdots .$

\noindent \textbf{Proof.}\ By multiplying in the Riordan group, we
immediately obtain
 \begin{equation*}
  \displaystyle E=D\ast P
  =\left( \frac{1}{1-x^{2}},\frac{x}{1-x^{2}}\right) \ast
   \left( \frac{1}{1+x},\frac{x}{1+x}\right)
  =\left( \frac{1}{1+x-x^{2}},\frac{x}{1+x-x^{2}}\right)
 \end{equation*}%
\noindent However
 \begin{equation*}
  \displaystyle[x^{n}]\frac{x^{k}}{(1+x-x^{2})^{k+1}}
  =(-1)^{n-k}\sum_{j=0}^{n-k}{\binom{k+j}{k}}{\binom{j}{n-j-k}}%
 \end{equation*}
Let
 $s_{n,k}=\sum_{j=0}^{n-k}{\binom{k+j}{k}}{\binom{j}{n-j-k}},$
we obtain the formula (5) and $s_{n,k}$ is the entry of the Riordan array
 \begin{equation*}
 \displaystyle E^{\ast}
 =\left( \frac{1}{1-x-x^{2}},\frac{x}{1-x-x^{2}}\right)
 .\quad _{\blacksquare }
 \end{equation*}
\noindent {\textbf{Remark.}}\ The matrices $E$ and $E^{\ast}$ could be called Fibonacci matrices.

\noindent {\textbf{Theorem~3.6.}}\ For $n\geq 0$, we have the
following formula
\begin{equation}
{m}_{n}=\sum_{k=0}^{n}{t}_{n,k}{s}_{k},
\end{equation}%
where the generating function of $\displaystyle({t}_{n,k})_{n\in
\mathbb{N}}$ is
 $\frac{1}{x}\left(\frac{\sqrt{1-2x+5x^{2}}+x-1}{2x}\right)^{k+1}, k=0,1,2,\cdots .$

\noindent \textbf{Proof.}\ The inverse of Riordan array $E$ is
 \begin{equation*}
  E^{-1}=\displaystyle
  \left(\frac{\sqrt{1-2x+5x^{2}}+x-1}{2x^{2}},\frac{\sqrt{1-2x+5x^{2}}+x-1}{2x}\right)
 \end{equation*}%
Giving the result. $\quad _{\blacksquare }$

\noindent {\textbf{4\hspace{0.4cm} Relations of the Generic Element
} }{\ }

In Section 3, we introduced two arrays $d_{n,k}$ and $s_{n,k}.$ In this section, we will discuss some identities related to them.
Let review them as follows firstly.

For the Riordan array
 $$D^{\ast}
 =\left( C\left(x^{2}\right) ,xC\left( x^{2}\right) \right)
 =\left[\begin{array}{ccccccc}
  {1} &  &  &  &  &  &    \\
  0 & {1} &  &  &  &  &    \\
  {1} & 0 & {1} &  &  &  &    \\
  0 & {2} & 0 & {1} &  &  &    \\
  {2} & 0 & {3} & {0} & 1 &  &    \\
  {0} & 5 & {0} & {4} & 0 & {1} &    \\
    &   & \cdots &  &  &  &  \ddots
 \end{array}\right], $$

\noindent $d_{n,k}$ is the generic element of $D^{\ast}$ and satisfies the
following recurrence relation
\begin{equation}
d_{n,k}=d_{n-1,k-1}+d_{n-1,k+1}.
\end{equation}
For the Riordan array
 $$E^{\ast}
 =\left( \frac{1}{1-x-x^{2}},\frac{x}{1-x-x^{2}}\right)
 =\left[\begin{array}{ccccccc}
  {1} &  &  &  &  &  &   \\
  1 & {1} &  &  &  &  &   \\
  {2} & 2 & {1} &  &  &  &  \\
  3 & {5} & 3 & {1} &  &  &  \\
  {5} & 10 & 9 & 4 & 1 &  &  \\
  8 & 20 & 22 & 14 & 5 & 1 & \\
   & & \cdots &  &  &  & \ddots
 \end{array}\right] ,$$

\noindent $s_{n,k}$ is the generic element of $E^{\ast }$ and
satisfies the following recurrence relation
\begin{equation}
s_{n,k}=s_{n-1,k-1}+s_{n-1,k}+s_{n-2,k}.
\end{equation}%

We need the combinatorial interpretation of entries in the above
matrices in terms of free symmetric lattice paths. To be precise,
a \emph{generalized free symmetric Motzkin path} of length $2n$ is a
path that the $i$-th step is up step (down step, $k$-horizontal
step) then the $2n-i$-th step is down step (up step,
$k$-horizontal step) without the restriction that it cannot go below
the $x$-axis. In this paper, when $k=0,1,2$, we call such a path a
free symmetric Dyck path, a free symmetric Motzkin path and a free
symmetric Schr\"{o}der path, respectively. Let $\mathcal{D}_{n},\mathcal{M}%
_{n},\mathcal{S}_{n}$, denote the set of free symmetric Dyck paths,
free symmetric Motzkin paths, free symmetric Schr\"{o}der paths of
length $2n$,
respectively. Note that $|\mathcal{D}_{n}|={2^{n}}$ while $|\mathcal{S}%
_{n}|=p_{n+1},$ the Pell numbers, (sequence A000129 in [8]) whose
generating function is $\frac{1}{1-2x-x^{2}}.$

A \emph{\ free symmetric MS path} is one that can have level steps
both of length one and two. Let ${\mathcal{MS}_{n}}$ denote the set
of  free symmetric MS paths  ending at $\left( 2n,0\right) . $
Moreover, ${|\mathcal{MS}_{n}|}=h_{n+1}$(sequence A006190 in [8])
whose generating function is $\frac{1}{1-3x-x^{2}}.$

Now we go back to the paths that don't go below the $x$-axis.

\noindent {\textbf{Theorem 4.1.}}\ Let $\mathcal{D}_{n,k}$ denote the
set of
symmetric Dyck paths of length $2n$ with the mid-height $k$. (The mid-height is the $y$-coordinate of the middle point.)  Then $d_{n,k}$ is the cardinality of $\mathcal{D}_{n,k}.$

\noindent \textbf{Proof.}\ Let $Q\in \mathcal{D}_{n,k}.$
If the $n$-th step of $Q$ is $U,$ then we can obtain a subpath $%
Q_{1}\in \mathcal{D}_{n-1,k-1}$ by deleting the $n$-th and $n+1$-th
steps of $Q.$ If the $n$-th step of $Q$ is $D,$ then we can obtain a subpath $%
Q_{2}\in \mathcal{D}_{n-1,k+1}$ by deleting the $n$-th and $n+1$-th
steps of $Q.$

So $|\mathcal{D}%
_{n,k}|=|\mathcal{D}_{n-1,k-1}|+|\mathcal{D}_{n-1,k+1}|.$ Combining
this with relation (7), we see that
$d_{n,k}=|\mathcal{D}_{n,k}|.\quad _{\blacksquare }$

\noindent {\textbf{Theorem 4.2.}}\ For $\displaystyle n\geq 0$ the sequence $%
{d}_{n,k}$ satisfies
\begin{equation}
\displaystyle\sum_{k=0}^{n}{d}_{n,k}={\binom{n}{\lfloor \frac{n}{2}\rfloor }}%
.
\end{equation}%
\noindent \textbf{Proof.} Consider the Riordan array $D^{\ast}$ and
the sequence $(1,1,1,\ldots )$ whose generating function is
$A(x)=\frac{1}{1-x},$ then
 \begin{equation*}
 \displaystyle\sum_{k=0}^{n}{d}_{n,k}
 =[x^{n}]C(x^{2})A\left( xC\left( x^{2}\right) \right)
 =[x^{n}]\frac{C(x^{2})}{1-xC\left(x^{2}\right)}
 ={\binom{n}{\lfloor \frac{n}{2}\rfloor }}%
 \end{equation*}

\noindent \textbf{Combinatorial proof.} Noting that the number of
symmetric
Dyck paths of length $2n$ equals $d_{n}={\binom{n}{\lfloor \frac{n}{2}%
\rfloor }},$ we have an immediate proof the theorem$.\quad
_{\blacksquare }$

\noindent {\textbf{Theorem 4.3.}}\ For $\displaystyle n\geq 0$ the sequence $%
{d}_{n,k}$ satisfies
\begin{equation}
\displaystyle\sum_{k=0}^{n}(k+1){d}_{n,k}=2^{n}.
\end{equation}%
\noindent \textbf{Proof.}\ Consider the Riordan array $D^{\ast }$
and
sequence $(1,2,3,\cdots )$ whose generating function is $A(x)=\frac{1}{%
(1-x)^{2}},$ then
 \begin{equation*}
 \begin{array}{ll}
 \displaystyle\sum_{k=0}^{n}(k+1){d}_{n,k}
 =[x^{n}]C(x^2)A(xC(x^2))
 =[x^{n}]\displaystyle\frac{C(x^{2})}{\left(1-xC(x^{2})\right)^{2}}
 =2^{n}\qquad
 \end{array}
 \end{equation*}

\noindent {\ \textbf{Combinatorial proof}.}\ To prove the theorem,
we will present a bijection between the set of restricted
$\mathcal{D}_{n,k}$ and the set $\mathcal{D}_{n}$. Clearly,
$|\mathcal{D}_{n}|=2^{n}.$ Suppose $Q\in \mathcal{D}_{n,k}.$ We will
present a set of maps $\phi _{i},i=0,1,2,\ldots k $ so that the left
hand side of identity (10) is the cardinality of the set
\begin{equation*}
\displaystyle\bigcup_{k=0}^{n}\bigcup_{Q\in \mathcal{D}_{n,k}}\{\phi
_{i}(Q):0\leq i\leq k\}.
\end{equation*}%
$\phi _{i}(Q)$ denotes the path obtained by changing each of the
last ascents to height $1,2,3,\ldots i$ to down steps. The
\emph{last ascent to
height $i$} of $Q$ is the last up step (going from the first step to the midpoint at the $n$%
-th step) which starts at height $i-1$ and ends at height $i$.

Since $Q$ has the mid-height $k$, it contains exactly $k$
last ascents of the left $n$ steps. By replacing the first $i$ last
ascents in $Q$ with $D$'s, we ensure that, to some point, the number
of $D$'s exceed the
number of preceding $U$'s, so that $\phi_{i}(Q)$ necessarily goes below the $%
x$-axis, eventually the $n$-th step ending at height $k-2i.$

See Figure 1 and 2 for an illustration of this map.

\setlength{\unitlength}{0.4cm}
\begin{picture}(30,8)
\put(2,0){\line(1,1){1}}\put(3,1){\line(1,-1){1}}\put(4,0){\line(1,1){2}}
\put(6,2){\line(1,-1){1}}\put(7,1){\line(1,1){2}}\put(9,3){\line(1,-1){1}}
\put(10,2){\line(1,1){2}}\put(12,4){\line(1,-1){2}}\put(14,2){\line(1,1){1}}
\put(15,3){\line(1,-1){1}}\put(16,2){\line(1,1){2}}
\put(18,4){\line(1,-1){2}}\put(20,2){\line(1,1){1}}\put(21,3){\line(1,-1){2}}
\put(23,1){\line(1,1){1}}\put(24,2){\line(1,-1){2}}\put(26,0){\line(1,1){1}}
\put(27,1){\line(1,-1){1}}
\put(3.0,1.0){\makebox(0,0)[c]{$\bullet$}}
\put(2.0,0.0){\makebox(0,0)[c]{$\bullet$}}\put(4.0,0.0){\makebox(0,0)[c]{$\bullet$}}
\put(5.0,1.0){\makebox(0,0)[c]{$\bullet$}}
\put(6.0,2.0){\makebox(0,0)[c]{$\bullet$}}\put(7.0,1.0){\makebox(0,0)[c]{$\bullet$}}
\put(8.0,2.0){\makebox(0,0)[c]{$\bullet$}}\put(9.0,3.0){\makebox(0,0)[c]{$\bullet$}}
\put(10.0,2.0){\makebox(0,0)[c]{$\bullet$}}\put(11.0,3.0){\makebox(0,0)[c]{$\bullet$}}
\put(12.0,4.0){\makebox(0,0)[c]{$\bullet$}}\put(13.0,3.0){\makebox(0,0)[c]{$\bullet$}}
\put(14.0,2.0){\makebox(0,0)[c]{$\bullet$}}\put(15.0,3.0){\makebox(0,0)[c]{$\bullet$}}
\put(16.0,2.0){\makebox(0,0)[c]{$\bullet$}}
\put(17.0,3.0){\makebox(0,0)[c]{$\bullet$}}\put(18.0,4.0){\makebox(0,0)[c]{$\bullet$}}
\put(19.0,3.0){\makebox(0,0)[c]{$\bullet$}}\put(20.0,2.0){\makebox(0,0)[c]{$\bullet$}}
\put(21.0,3.0){\makebox(0,0)[c]{$\bullet$}}
\put(22.0,2.0){\makebox(0,0)[c]{$\bullet$}}\put(23.0,1.0){\makebox(0,0)[c]{$\bullet$}}
\put(24.0,2.0){\makebox(0,0)[c]{$\bullet$}}\put(25.0,1.0){\makebox(0,0)[c]{$\bullet$}}
\put(26.0,0.0){\makebox(0,0)[c]{$\bullet$}}\put(27.0,1.0){\makebox(0,0)[c]{$\bullet$}}
\put(28.0,0.0){\makebox(0,0)[c]{$\bullet$}}

\put(3.9,0.8){\makebox(0,0)[l]{$\ast$}}\put(6.9,1.8){\makebox(0,0)[l]{$\ast$}}
\put(13.9,2.8){\makebox(0,0)[l]{$\ast$}}

\end{picture}\newline

\noindent Figure 1: $Q=U\ D\ \mathbf{U}\ U\ D\ \mathbf{U}\ U\ D\ U\
U\ D\ D\ \mathbf{U}\ D\ U\ U\ D\ D\ U\ D\ D\ U\ D\ D\ U\ D,$ a
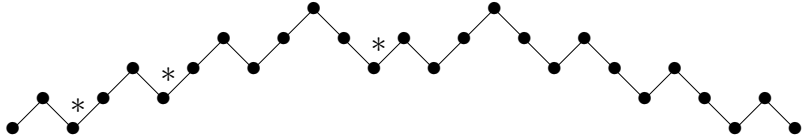
restricted symmetric Dyck path of length $2n=26$ and the mid-height $k=3.$ Last ascents are indicated with a bold $\mathbf{U}$
and marked above with an single star.

\setlength{\unitlength}{0.4cm}
\begin{picture}(14,4)
\put(0,0){\line(1,1){1}}\put(1,1){\line(1,-1){2}}\put(3,-1){\line(1,1){1}}
\put(4,0){\line(1,-1){2}}\put(6,-2){\line(1,1){1}}\put(7,-1){\line(1,-1){1}}
\put(8,-2){\line(1,1){2}}\put(10,0){\line(1,-1){2}}\put(12,-2){\line(1,1){1}}
\put(13,-1){\line(1,-1){1}}\put(14,-2){\line(1,1){2}}\put(16,0){\line(1,-1){2}}
\put(18,-2){\line(1,1){1}}\put(19,-1){\line(1,-1){1}}\put(20,-2){\line(1,1){2}}
\put(22,0){\line(1,-1){1}}\put(23,-1){\line(1,1){2}}\put(25,1){\line(1,-1){1}}

\put(1.0,1.0){\makebox(0,0)[c]{$\bullet$}}
\put(0.0,0.0){\makebox(0,0)[c]{$\bullet$}}\put(2.0,0.0){\makebox(0,0)[c]{$\bullet$}}\put(3.0,-1.0){\makebox(0,0)[c]{$\bullet$}}
\put(4.0,0.0){\makebox(0,0)[c]{$\bullet$}}\put(5.0,-1.0){\makebox(0,0)[c]{$\bullet$}}
\put(6.0,-2.0){\makebox(0,0)[c]{$\bullet$}}\put(7.0,-1.0){\makebox(0,0)[c]{$\bullet$}}
\put(8.0,-2.0){\makebox(0,0)[c]{$\bullet$}}\put(9.0,-1.0){\makebox(0,0)[c]{$\bullet$}}
\put(10.0,0.0){\makebox(0,0)[c]{$\bullet$}}\put(11.0,-1.0){\makebox(0,0)[c]{$\bullet$}}
\put(12.0,-2.0){\makebox(0,0)[c]{$\bullet$}}\put(13.0,-1.0){\makebox(0,0)[c]{$\bullet$}}

\put(14.0,-2.0){\makebox(0,0)[c]{$\bullet$}}
\put(15.0,-1.0){\makebox(0,0)[c]{$\bullet$}}\put(16.0,0.0){\makebox(0,0)[c]{$\bullet$}}\put(17.0,-1.0){\makebox(0,0)[c]{$\bullet$}}
\put(18.0,-2.0){\makebox(0,0)[c]{$\bullet$}}\put(19.0,-1.0){\makebox(0,0)[c]{$\bullet$}}
\put(20.0,-2.0){\makebox(0,0)[c]{$\bullet$}}\put(21.0,-1.0){\makebox(0,0)[c]{$\bullet$}}
\put(22.0,0.0){\makebox(0,0)[c]{$\bullet$}}\put(23.0,-1.0){\makebox(0,0)[c]{$\bullet$}}
\put(24.0,0.0){\makebox(0,0)[c]{$\bullet$}}\put(25.0,1.0){\makebox(0,0)[c]{$\bullet$}}
\put(26.0,0.0){\makebox(0,0)[c]{$\bullet$}}

\put(1.9,-0.8){\makebox(0,0)[l]{$\ast$}}\put(4.9,-1.8){\makebox(0,0)[l]{$\ast$}}
\end{picture}\newline

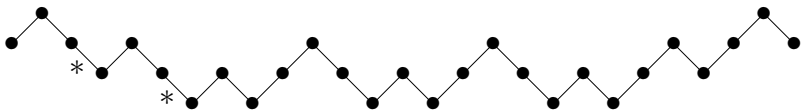
$\newline $ \noindent Figure 2: $\phi _{2}(Q)=U\ D\ \mathbf{D}\ U\
D\ \mathbf{D}\ U\ D\ U\ U\ D\ D\ U,$ the unrestricted Dyck path of
length $n=13$ and terminal height $j=k-2i=3-2(2)=-1$ which is
obtained by changing the last ascents 1 and 2 to down steps. Premier
descents are indicated with a bold $\mathbf{D}$ and marked with a
single star .

From the construction, we see directly that $\phi _{i}$ is one to
one. In order to show that the set $\{Q,\phi _{0}(Q),\ldots \phi
_{k}(Q):Q\in
\mathcal{D}_{n,k},0\leq k\leq n\}$ is indeed equal to the set $\mathcal{D}%
_{n},$ we need only show that it is possible to take a path $p\in \mathcal{D}%
_{n}$ and recover its unique preimage $Q$ under $\phi _{i}$ for some
$i.$

Let $p\in \mathcal{D}_{n}$ with mid-height $j$. We
consider the left $n$ steps. Since $p$ goes below the $x$-axis there exists a first step down from height $m$ to $m-1$ from $m=0,-1,-2,\ldots$ as one proceeds from left to right. We call these premier steps. Assuming that  the number of the first step below
the $x$-axis of $p$ is $i(>0)$ , then changing each of the $i$
premier steps to up steps will create a new path $Q$ which stays
above the $x$-axis and ends at height $k=j+2i.$

So this map is a bijection. $\quad _{\blacksquare }$

\noindent {\ \textbf{Corollary 4.4}.}\ The average mid-height of the symmetric Dyck paths from $%
\left( 0,0\right) $ to $\left( 2n,0\right) $  is%
\begin{equation*}
\left\{
\begin{array}{cc}
\frac{2^{2m}-\binom{2m}{m}}{\binom{2m}{m}}\sim \sqrt{\pi m}-1\sim \sqrt{\pi m%
} & \text{ if }n=2m \\
\frac{2^{2m+1}-\binom{2m+1}{m}}{\binom{2m+1}{m}}\sim \sqrt{\pi
m}-1\sim
\sqrt{\pi m} & \text{if }n=2m+1%
\end{array}%
\right.
\end{equation*}

\noindent {\ \textbf{Proof}.}\ To compute the total height we multiply the matrix
$D^{\ast}$ by the column vector
 $\left( 0,1,2,3,\cdots \right)^{T}
 =\left(1,2,3,\cdots\right) ^{T}-\left( 1,1,1,\cdots \right) ^{T}.$
We consider the case $n=2m$, the other being quite similar. Combining the results of Theorem 4.2 and 4.3 and
 $n!\sim \left(\frac{n}{e}\right)^n\sqrt{2n\pi}(n\rightarrow+\infty)$
gives the result. $\quad_ \blacksquare $

Thus the average mid-height tends to infinity as $n$
gets large. How about the number of points where the path touches
the $x$-axis?

\noindent\ \textbf{Theorem 4.5.} The generating function for the total
number of points on the $x$-axis for paths from $\left( 0,0\right) $
to $\left( 2n,0\right) $ is
\begin{eqnarray*}
&&2\left[ C\left( x^{2}\right) \cdot \frac{1}{2x}\left( \sqrt{\frac{1+2x}{%
1-2x}}-1\right) \right] -C\left( x^{2}\right)  \\
&=&1+2x+5x^{2}+8x^{3}+18x^{4}+30x^5+65x^{6}+112z^{7}+\cdots .
\end{eqnarray*}%
\noindent For $n=2m$ the generating function of total number of
points on the $x-$axis is
\begin{eqnarray*}
&C\left(
x^{2}\right)[2B(x^2)-1]=&1+5x^{2}+18x^{4}+65x^{6}+238x^{8}+\cdots .
\end{eqnarray*}%
and thus the average number of points on the $x-$axis is
\begin{equation*}
\begin{array}{c}
\frac{\binom{2m+2}{m+1}-\frac{1}{m+1}\binom{2m}{m}}{\binom{2m}{m}}=\frac{4m+1%
}{m+1}\rightarrow 4%
\end{array}%
\end{equation*}%
while if $n=2m+1$  then the corresponding results are
\begin{equation*}
2xB\left( x^{2}\right) C\left( x^{2}\right) ^{2}=
2x+8x^{3}+{30x^{5}}+112x^{7}+\cdots
\end{equation*}%
and the average number of points is
\begin{equation*}
\frac{2\binom{2m+2%
}{m}}{\binom{2m+1}{n}}=\frac{2\binom{2m+2%
}{m}}{\frac{1}{2}\binom{2m+2}{m}}=4\cdot \frac{n+1}{n+2}\rightarrow
4.\quad \blacksquare
\end{equation*}%

\noindent{\textbf{Remark 1.}} When $n$ is even the total number of points on the
$x$-axis is $\left( 4m+1\right) C_{m}.$

\noindent{\textbf{Remark 2.}} Without the symmetric condition the classical theorem  of Dershowitz and Zaks [2], stated for ordered trees, gives the total number of points on the $x$-axis is $C_{n+1}$ and thus the average value becomes $\frac{4n+2}{n+2}\rightarrow 4$

\noindent{\textbf{Theorem 4.6.}} Let $\mathcal{S}_{n,k}$ denote the
set of all free symmetric Schr\"{o}der paths of length $2n$ that the
left $n$ steps
contain $k$ up steps. Then $s_{n,k}$ enumerates the cardinality of $\mathcal{%
S}_{n,k}.$

\noindent\textbf{Proof.}\ Let $Q\in \mathcal{S}_{n,k}.$ If the middle step is $U,$
 then we can obtain a subpath $Q_{1}\in \mathcal{S}%
_{n-1,k-1}$ by deleting the middle step and the following step of
$Q$. If the middle step is $D,$ then we can obtain a subpath $Q_{2}\in \mathcal{S}%
_{n-1,k}$ by deleting the middle step and the following step of $Q$. If the middle step
is $H,$ then we can obtain a subpath $Q_{3}\in \mathcal{S}%
_{n-2,k}$ by deleting the middle step and the following step of $Q$.

So $|\mathcal{S}%
_{n,k}|=|\mathcal{S}_{n-1,k-1}|+|\mathcal{S}_{n-1,k}|+|\mathcal{S}_{n-2,k}|.$
Combined with the relation (8), we can see that
$s_{n,k}=|\mathcal{S}_{n,k}|.\quad _{\blacksquare}$

\noindent{\textbf{Theorem 4.7.}}\ For $\displaystyle n\geq 0$ the sequence ${%
s}_{n,k}$ satisfies
\begin{align}
\displaystyle \sum_{k=0}^{n}{s}_{n,k}= p_{n+1}.
\end{align}
\noindent{\ \textbf{Proof}.} \noindent By FTRA and $A(x)=\frac{1}{1-x}$  we have
\begin{equation*}
  \sum_{k=0}^{n}{s}_{n,k}=[x^{n}]
  \left(\frac {1}{1-x-x^{2}},\frac
  {x}{1-x-x^{2}}\right)*\frac{1}{1-x}
 =[x^{n}]\frac {1}{1-2x-x^{2}}
 =p_{n+1}
\end{equation*}
\noindent This completes the proof.

\noindent{\ \textbf{Combinatorial proof}.}\ Clearly, the term $\sum_{k=0}^{n}%
{s}_{n,k}$ suggests the number of the free symmetric Schr\"{o}der
paths of length $2n$. This observation gives an immediate proof of the identity.  $\quad_{\blacksquare}$

\noindent{\textbf{Theorem 4.8.}}\ For $\displaystyle n\geq 0$ the sequence ${%
s}_{n,k}$ satisfies
\begin{align}
\displaystyle\sum_{k=0}^{n}2^{k}s_{n,k}=h_{n+1}.
\end{align}
\noindent{\ \textbf{Proof}.} By FTRA and $A(x)= \frac{1}{1-2x}$ we have
\begin{equation*}
  \displaystyle\sum_{k=0}^{n}2^{k}s_{n,k}
 =[x^{n}]\left(\frac {1}{1-x-x^{2}},\frac
  {x}{1-x-x^{2}}\right)*\frac{1}{1-2x}
 =[x^{n}]\frac {1}{1-3x-x^{2}}
 =h_{n+1}
\end{equation*}
\noindent This completes the proof.

\noindent{\ \textbf{Combinatorial proof}.}\ Suppose $Q\in \mathcal{S}%
_{n,k}(k>0).$ We will present a set of maps $\phi_{j_{1},j_{2}\ldots
j_{i}\ldots j_{k}}(j_{i}=0$ or 1 ) so that the left hand side of
(12) is the cardinality of the set
\begin{equation*}
\displaystyle \mathcal{S}_{n,0}\bigcup
\bigcup_{k=1}^{n}\bigcup_{Q\in
\mathcal{S}_{n,k}}\{\phi_{j_{1},j_{2}\ldots j_{i}\ldots j_{k}(Q)}:
j_{i}=0 \ \mathrm{or} \ 1\} .
\end{equation*}
$\phi_{j_{1},j_{2}\ldots j_{i}\ldots j_{k}(Q)}$ denote the path
obtained by changing each of the $k$ up steps to horizontal steps
$h$ or remain $U$. If we change the up step $i$ to horizontal step,
we have $j_{i}=1.$ If the up step $i$ remains we have $j_{i}=0.$
See the following for an illustration of this map.

\setlength{\unitlength}{0.5cm}
\begin{picture}(30,3)
\put(3,0){\line(1,1){1}}\put(4,1){\line(1,-1){1}}\put(5,0){\line(1,-1){1}}
\put(6,-1){\line(1,1){1}}\put(7,0){\line(1,0){2}}\put(9,0){\line(1,1){1}}
\put(10,1){\line(1,-1){1}}\put(11,0){\line(1,1){1}}\put(12,1){\line(1,-1){1}}
\put(13,0){\line(1,0){2}}\put(15,0){\line(1,-1){1}}
\put(16,-1){\line(1,1){2}}\put(18,1){\line(1,-1){1}}
\put(4.0,1.0){\makebox(0,0)[c]{$\bullet$}}
\put(3.0,0.0){\makebox(0,0)[c]{$\bullet$}}\put(5.0,0.0){\makebox(0,0)[c]{$\bullet$}}\put(6.0,-1.0){\makebox(0,0)[c]{$\bullet$}}
\put(7.0,0.0){\makebox(0,0)[c]{$\bullet$}}
\put(9.0,0.0){\makebox(0,0)[c]{$\bullet$}}\put(10.0,1.0){\makebox(0,0)[c]{$\bullet$}}
\put(11.0,0.0){\makebox(0,0)[c]{$\bullet$}}\put(12.0,1.0){\makebox(0,0)[c]{$\bullet$}}
\put(13.0,0.0){\makebox(0,0)[c]{$\bullet$}}
\put(15.0,0.0){\makebox(0,0)[c]{$\bullet$}}\put(16.0,-1.0){\makebox(0,0)[c]{$\bullet$}}\put(17.0,0.0){\makebox(0,0)[c]{$\bullet$}}
\put(18.0,1.0){\makebox(0,0)[c]{$\bullet$}}\put(19.0,0.0){\makebox(0,0)[c]{$\bullet$}}

\put(3.1,0.8){\makebox(0,0)[l]{$\ast$}}\put(6.1,-0.3){\makebox(0,0)[l]{$\ast$}}
\put(9.1,0.8){\makebox(0,0)[l]{$\ast$}}

\end{picture}\newline

\noindent Figure 1 : $Q=\mathbf{U}\ D\ D\ \mathbf{U}\ H\ \mathbf{U}\
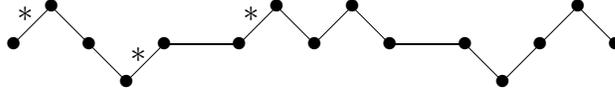
D\ U\ D\ H\ D\ U\ U\ D,$ a free symmetric Schr\"{o}der path of
length $2n=16$ such that the left $n$ steps contain $k=3$ up steps. The
$k$ up steps are indicated with a bold $\mathbf{U}$ and marked above
with an star.

\setlength{\unitlength}{0.5cm}
\begin{picture}(30,2)
\put(3,0){\line(1,0){1}}\put(4,0){\line(1,-1){2}}
\put(6,-2){\line(1,1){1}}\put(7,-1){\line(1,0){2}}\put(9,-1){\line(1,0){1}}
\put(10,-1){\line(1,-1){1}}\put(11,-2){\line(1,1){1}}\put(12,-1){\line(1,0){1}}
\put(13,-1){\line(1,0){2}}\put(15,-1){\line(1,-1){1}}
\put(16,-2){\line(1,1){2}}\put(18,0){\line(1,0){1}}
\put(4.0,0.0){\makebox(0,0)[c]{$\bullet$}}
\put(3.0,0.0){\makebox(0,0)[c]{$\bullet$}}\put(5.0,-1.0){\makebox(0,0)[c]{$\bullet$}}
\put(6.0,-2.0){\makebox(0,0)[c]{$\bullet$}}
\put(7.0,-1.0){\makebox(0,0)[c]{$\bullet$}}
\put(9.0,-1.0){\makebox(0,0)[c]{$\bullet$}}\put(10.0,-1.0){\makebox(0,0)[c]{$\bullet$}}
\put(11.0,-2.0){\makebox(0,0)[c]{$\bullet$}}\put(12.0,-1.0){\makebox(0,0)[c]{$\bullet$}}
\put(13.0,-1.0){\makebox(0,0)[c]{$\bullet$}}
\put(15.0,-1.0){\makebox(0,0)[c]{$\bullet$}}\put(16.0,-2.0){\makebox(0,0)[c]{$\bullet$}}
\put(17.0,-1.0){\makebox(0,0)[c]{$\bullet$}}
\put(18.0,0.0){\makebox(0,0)[c]{$\bullet$}}\put(19.0,0.0){\makebox(0,0)[c]{$\bullet$}}

\put(3.1,0.3){\makebox(0,0)[l]{$\ast$}}
\put(9.2,-0.5){\makebox(0,0)[l]{$\ast$}}

\end{picture}\newline

\quad

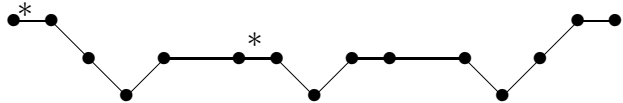
\noindent Figure 2: $\phi_{1,0,1}(Q)=\mathbf{h}\ D\ D\ {U}\ H\
\mathbf{h}\
D\ U\ h\ H\ D\ U\ U\ h,$ the free generalized Schr\"{o}der path of length $%
2n=16$ which obtained by changing the up steps 1 and 3 to horizontal steps $%
h.$ The changed horizontal steps are indicated with a bold $h$ and
marked above with an star.

Let $P\in  {\mathcal{MS}_{n}}$. We concentrate on the left half of
the path which has $j$ horizontal steps and $i$ up steps. Then
changing each of
the $j$ horizontal steps $h$ to up steps will get a free symmetric Schr\"{o}%
der path of length $2n$ such that the left $n$ steps contain $k=i+j$ up
steps.

From above, the set $\displaystyle \mathcal{S}_{n,0}\bigcup
\bigcup_{k=1}^{n}\bigcup_{Q\in
\mathcal{S}_{n,k}}\{\phi_{j_{1},j_{2}\ldots
j_{i}\ldots j_{k}}: j_{i}=0 \ \mathrm{or} \ 1\}$ is indeed equal to the set $%
{\mathcal{MS}_{n}}\quad _{\blacksquare}$

{\small {\ }}

{\small \ }

\end{document}